\documentclass[11pt]{article}
\PassOptionsToPackage{obeyspaces}{url}
\usepackage[hidelinks]{hyperref}
\usepackage{bookmark}
\bookmarksetup{open,numbered,depth=5}
\usepackage{amsfonts}
\usepackage{amsmath}
\usepackage{amssymb, amsthm, enumitem}
\usepackage{breakurl}
\usepackage{pgf,tikz}
\usepackage{bm}
\usepackage{cancel}
\usepackage{comment}
\usetikzlibrary{calc}

\usepackage{etoolbox} 
\patchcmd{\thebibliography}{\leftmargin\labelwidth}{\leftmargin\labelwidth\addtolength\itemsep{-0.1\baselineskip}}{}{}

\usepackage{mathtools}
\usepackage{xstring}

\author{
Ting-Wei Chao\thanks{Department of Mathematics, Massachusetts Institute of Technology, Cambridge, MA, USA. \texttt{twchao@mit.edu}} \and Zichao Dong\thanks{Extremal Combinatorics and Probability Group (ECOPRO), Institute for Basic Science (IBS), Daejeon, South Korea. Supported by the Institute for Basic Science (IBS-R029-C4). \texttt{zichao@ibs.re.kr}. } 
\and Zijun Shen\thanks{School of Mathematical Sciences, Peking University, Beijing 100871. \texttt{2300010734@stu.pku.edu.cn}. }
\and Ningyuan Yang\thanks{School of Mathematical Sciences, Fudan University, Shanghai 200433. \texttt{nyyang23@m.fudan.edu.cn}}}

\title{Many cliques with small degree powers}
\date{}

\usepackage[nameinlink]{cleveref}

\newtheorem{theorem}{Theorem}
\newtheorem{lemma}[theorem]{Lemma}

\newtheorem{proposition}[theorem]{Proposition}

\newtheorem{claim}{Claim}

\newtheorem{question}{Question}

\newcommand*{\eqdef}{\stackrel{\mbox{\normalfont\tiny def}}{=}} 

\newcommand*{\p}{\mathbb{P}}
\newcommand*{\N}{\mathbb{N}}                                    
\newcommand*{\Z}{\mathbb{Z}}                                    
\newcommand*{\R}{\mathbb{R}}                                    

\newcommand*{\bH}{\mathbb{H}}
\newcommand*{\bP}{\mathbb{P}}
\newcommand\coin{
\mathrel{\ooalign{\hss$\bigcirc$\hss\cr\kern0.7ex\hbox{\scalebox{0.8}{$\$$}}}}\,}

\DeclareMathOperator{\Bern}{Ber}
\DeclareMathOperator{\supp}{supp}

\allowdisplaybreaks[4]

\begin{document}

\maketitle

\begin{abstract}
    Suppose $0 < p \le \infty$. For a simple graph $G$ with a vertex-degree sequence $d_1, \dots, d_n$ satisfying $(d_1^p + \dots + d_n^p)^{1/p} \le C$, we prove asymptotically sharp upper bounds on the number of $t$-cliques in $G$. This result bridges the $p = 1$ case, which is the notable Kruskal--Katona theorem, and the $p = \infty$ case, known as the Gan--Loh--Sudakov conjecture, and resolved by Chase. In particular, we demonstrate that the extremal construction exhibits a dichotomy between a single clique and multiple cliques at $p_0 = t - 1$. Our proof employs the entropy method. 
\end{abstract}

\section{Introduction} \label{sec:intro}

\paragraph{Background.} In extremal combinatorics, an archetype problem is to study the maximum number of subgraphs isomorphic to a fixed graph under certain restrictions. A classical result in this line is the Kruskal--Katona theorem \cite{kruskal,katona}, which is of central significance in extremal set theory. Over the years, the Kruskal--Katona theorem has also shown many connections to related fields such as extremal graph theory and discrete geometry. See, for example, \cite{frankl_tokushige,chao_yu,chao2023tightboundstructuraltheorem}. 

\smallskip

Throughout this paper, every graph or hypergraph is simple. 

Let $e$ be the number of edges and $\Delta$ be the maximum degree of a graph $G$. We consider the following extremal questions concerning $k_3(G)$, the number of triangle subgraphs in $G$: 

\begin{question} \label{ques:triangle_sumdeg}
    For fixed $e$, what is the maximum possible value of $k_3(G)$? 
\end{question}

\begin{question} \label{ques:triangle_maxdeg}
    For fixed $\Delta$, what is the maximum possible value of $k_3(G)$? 
\end{question}

For \Cref{ques:triangle_sumdeg}, the Kruskal--Katona theorem tells us that $k_3(G)$ is maximized roughly when $G$ is a large clique. In this paper, we focus on understanding the asymptotical behavior of $k_3(G)$, and so by ``Kruskal--Katona theorem'' we usually refer to the corollary below of its approximate version, which was formulated by Lov\'{a}sz. See \cite[Chapter 13, Exercise 31b]{lovasz} for a reference. 

\begin{theorem} \label{thm:kruskalkatona_lovasz}
    Let $t \ge 3$ be an integer, and suppose $u \ge t$ is a real number. Then every graph with exactly $\binom{u}{2}$ edges contains at most $\binom{u}{t}$ cliques of size $t$. 
\end{theorem}

In particular, upon setting $t = 3$ \Cref{thm:kruskalkatona_lovasz} answers \Cref{ques:triangle_sumdeg}. 

\smallskip

For \Cref{ques:triangle_maxdeg}, one has to impose a restriction that $G$ is an $n$-vertex graph, for otherwise the answer is obviously infinite. Chase \cite{chase} proved that $k_3(G)$ is maximized roughly when $G$ is a disjoint union of $(\Delta + 1)$-cliques. In fact, Chase \cite{chase} established the following result: 

\begin{theorem} \label{thm:GLS_chase}
    For any positive integers $n, \Delta$, any graph on $n$ vertices of maximum degree $\Delta$ contains at most $q\binom{\Delta+1}{t} + \binom{r}{t}$ cliques of size $t$, where $n = q(\Delta+1) + r$ for integers $q$ and $r \in [0, \Delta]$. 
\end{theorem}

Historically, Galvin \cite{galvin} initiated the study of maximum possible number of independent sets in an $n$-vertex graph $G$ of minimum degree $\delta$ (see also \cite{cutler_radcliffe}). Then Engbers and Galvin asked what happens if we would like to maximize the number of independent sets of some fixed size, given the minimum degree of the graph. By taking the complement graph, Gan, Loh, and Sudakov \cite{gan_loh_sudakov} conjectured that $k_t(G)$, the number of $t$-cliques in $G$, is maximized when $G$ is a disjoint union of $(\Delta+1)$-cliques and a smaller clique, given that the maximum degree of $G$ is $\Delta$. Indeed, they managed to prove using the Kruskal--Katona theorem that it suffices to resolve the conjecture for $t = 3$. A couple years later, the Gan--Loh--Sudakov conjecture was resolved by Chase \cite{chase} via a nice elementary approach. Later on, Chao and Dong \cite{chao_dong} generalized Chase's method and gave a unified proof of the Gan--Loh--Sudakov conjecture without applying the Kruskal--Katona theorem. 

\paragraph{The generalization.} For a graph $G$ on vertices $v_1, \dots, v_n$, we call ${\bm d}(G) \eqdef \bigl( \deg(v_1), \dots, \deg(v_n) \bigr)$ the degree sequence of $G$. Suppose $p \in (0, \infty)$ and refer to the \emph{degree sequence $\ell^p$-norm} of $G$ as
\[
\|\bm{d}(G)\|_p \eqdef \biggl( \sum_{i=1}^n \deg(v_i)^p \biggr)^{1/p}. 
\]
Also write $\|\bm{d}(G)\|_{\infty} \eqdef \lim\limits_{p \to \infty} \|\bm{d}(G)\|_p$, which is equal to the maximum degree in graph $G$. In fact, the terminology is abused here, since this quantity does not satisfy the norm axioms when $0 < p < 1$. 

Within the notations above, \Cref{ques:triangle_sumdeg} and \Cref{ques:triangle_maxdeg} ask the maximum value of $k_3(G)$, given that $\|\bm{d}(G)\|_1$ or $\|\bm{d}(G)\|_{\infty}$ is fixed, respectively. It is natural to interpolate between these two versions and ask the same question for intermediate $\ell^p$-norms (including the $0 < p < 1$ case): 

\begin{question} \label{ques:triangle_pnorm}
    Let $G$ be a graph with fixed $\|\bm{d}(G)\|_p$. What is the maximum of $k_3(G)$? 
\end{question}

In this paper, we carry out a systematic study on \Cref{ques:triangle_pnorm}, concerning not only the number of triangles but general $t$-cliques as well. In the literature, Tur\'{a}n-type problems ask the maximum number of edges (equivalently, the maximum of $\|\bm{d}(G)\|_1$) when an $n$-vertex graph contains no given subgraphs. Caro and Yuster \cite{caro_yuster} initiated the study of degree power Tur\'{a}n-type problems, where the $\|\bm{d}(G)\|_1$ norm is replaced by $\|\bm{d}(G)\|_p$. See \cite{bollobas_nikiforov_2004,pikhurko_taraz,nikiforov,bollobas_nikiforov_2012,gu_li_shi,gerbner} for various results in this line. 

\paragraph{Main result.} Recall that the extremal graph for \Cref{ques:triangle_sumdeg} is one large clique, while the extremal graph for \Cref{ques:triangle_maxdeg} is a union of small cliques. So, intuitively one can expect a threshold $p_0$ such that the extremal graph for \Cref{ques:triangle_pnorm} varies from a single clique to a disjoint union of cliques as $p$ goes above $p_0$. Our main result, \Cref{thm:main}, is an asymptotically tight upper bound on the number of $t$-cliques (not only the number of triangles) in a graph of fixed degree sequence $\ell^p$-norm for every $p \ge 1$. As we shall see in the main theorem, the threshold is exactly $p_0 = t - 1$. 

Before stating the theorem, we need to introduce some definitions. Define the function
\[
h(x) = \frac{1}{x(x-1)^p} \binom{x}{t}
\]
on the interval $(t-1, +\infty)$. When $p > t - 1$, set $s_{\R}$ to be the real number maximizing $h(x)$, and set $s_{\N}$ to be an arbitrary integer maximizing $h(x)$ among the positive integers $\{t, t+1, t+2, \dots\}$. We shall see (in \Cref{sec:plarge}) that $s_{\R}$ is unique (well-defined), and there are at most two candidates for $s_{\N}$. 
        
\begin{theorem} \label{thm:main}
    Suppose $C > 0$. Let $G$ be a graph with $\|\bm{d}(G)\|_p \le C$. 
    \vspace{-0.5em}
    \begin{enumerate}[label=(\roman*), ref=(\roman*)]
        \item \label{main_psmall} If $0 < p \le t - 1$, then $k_t(G) \le \max \bigl\{ \binom{u}{t}, 0 \bigr\}$, where $u$ is the positive real with $C = \sqrt[p]{u}(u-1)$.
        \vspace{-0.5em}
        \item \label{main_plarge} If $p > t - 1$, then $k_t(G) \le \frac{C^p}{u(u-1)^p} \binom{u}{t}$, where $u=s_{\N}$.
    \end{enumerate}
    \vspace{-0.5em}
\end{theorem}

It is not hard to see that the $u$ satisfying $C = \sqrt[p]{u}(u-1)$ is unique in the first part. 

We will also discuss the regime when the number of vertices in $G$ is fixed while $p > t - 1$. In the case $0 < p \le t-1$, this question has the same answer as in  \Cref{thm:main}\ref{main_psmall}. 

\begin{theorem} \label{thm:clique_plarge}
    Suppose $p > t - 1$ and $n(s_{\R}-1)^p \le C^p$. If $G$ is an $n$-vertex graph with $\|\bm{d}(G)\|_p \le C$, then $k_t(G) \le \frac{n}{u}{\binom{u}{t}}$, where $u = \frac{C}{\sqrt[p]{n}} + 1$. 
\end{theorem}

In fact, \Cref{thm:main,thm:clique_plarge} are best possible under certain number theoretical restrictions: 
\vspace{-0.5em}
\begin{itemize}
    \item In \Cref{thm:main}\ref{main_psmall}, the sharpness is witnessed by a $u$-clique, provided that $u \in \Z$, and
    \vspace{-0.5em}
    \item in \Cref{thm:main}\ref{main_plarge}, when $\frac{C^p}{u(u-1)^p} \in \Z$ that is witnessed by $\frac{C^p}{u(u-1)^p}$ many disjoint $u$-cliques.
    \vspace{-0.5em}
    \item The tightness of \Cref{thm:clique_plarge}, given that $u, \frac{n}{u} \in \Z$, is witnessed by $\frac{n}{u}$ many disjoint $u$-cliques. 
\end{itemize}
\vspace{-0.5em}

\smallskip

Balogh, Clemen, and Lidick\'{y} \cite{balogh_clemen_lidicky_2022a,balogh_clemen_lidicky_2022b} recently brought up the study of hypergraph degree power Tur\'{a}n-type problems, with a highlight on asymptotically determining the maximum of $\|\bm{d}(G)\|_2$ in tetrahedron-free graphs. Chen et al.~\cite{chen_ilkovic_leon_liu_pikhurko} extensively studied the $(t, p)$-norm Tur\'an-type problems concerning hypergraphs. We shall also extend our results to hypergraphs (\Cref{thm:main_hyper}). 

\paragraph{Proof overview.} Upon setting $p = 1$, \Cref{thm:main}\ref{main_psmall} is equivalent to \Cref{thm:kruskalkatona_lovasz}. So, the first part itself can be viewed as a generalization of the notable Kruskal--Katona theorem. However, it seems difficult to deduce \Cref{thm:main}\ref{main_psmall} directly from \Cref{thm:kruskalkatona_lovasz}. Also, the classical combinatorial shifting method due to Frankl \cite{frankl} does not work here since the function $f(x) = x^p$ is strictly convex when $p > 1$ (unlike the $p = 1$ boundary case). To overcome the aforementioned difficulties, our approach is to employ entropy estimates which were recently applied to related problems by Chao and Yu \cite{chao_yu}. We also use this method to prove \Cref{thm:clique_plarge}. \Cref{thm:main}\ref{main_plarge} follows from a counting argument.

\paragraph{Paper organization.} We set off by introducing the necessary definitions and useful properties on entropy in \Cref{sec:entropy}. We prove \Cref{thm:main}\ref{main_psmall} via entropy estimates in \Cref{sec:psmall}. In \Cref{sec:plarge}, we first show that the critical points $s_{\R}, s_{\N}$ are well-defined by analyzing the function $h(x)$, and present the proofs of \Cref{thm:clique_plarge} and \Cref{thm:main}\ref{main_plarge}. Finally, we discuss further generalizations of our results to hypergraphs in \Cref{sec:remark} and conclude with some potential directions for future work.

\section{Entropy preliminaries}\label{sec:entropy}

As mentioned in the introduction, we shall apply the entropy method in the proofs of \Cref{thm:main}\ref{main_psmall} and \Cref{thm:clique_plarge}. We begin with some related terminologies and basic results. 

In order to prove \Cref{thm:main}\ref{main_psmall}, we collect some necessary definitions and properties of entropy. For an exposition of entropy in combinatorics, see \cite[Section~15.7]{alon_spencer}. For more recent applications of entropy applied to subgraph counting and Kruskal--Katona type problems, see \cite{chao2023tightboundstructuraltheorem,chao2024purelyentropicapproachrainbow}. 

Throughout this paper, every probability space $(\Omega, \p)$ is \emph{discrete} (i.e., $|\Omega| < \infty$). The support of a discrete random variable $X$ is defined to be the set of all $x\in X(\Omega)=\bigl\{ X(\omega) : \omega \in \Omega \bigr\}$ such that the probability $p_X(x)\eqdef\p(X=x)>0$, and we denote by $\supp(X)$ the support of $X$. Every random variable has \emph{finite support} (i.e.,~$\lvert \supp(X) \rvert < \infty$). The \emph{entropy} (also known as the \emph{Shannon entropy}) of a random variable $X$ is defined as
\[
\bH(X) \eqdef \sum_{x \in \supp(X)} -p_X(x) \log_2 p_X(x). 
\]
Informally speaking, the entropy measures in bits the amount of information carried by an random variable. For instance, the entropy of a random bit (i.e., a fair coin toss) is $1$. Denote by $\Bern(p)$ the Bernoulli random variable of parameter $p$. Then $\bH\bigl(\Bern(p)\bigr) = -p \log_2 p - (1-p) \log_2 (1-p)$. 

\smallskip

The following uniform bound of entropy is important. 

\begin{proposition} \label{prop:entropy_max}
    Let $X$ be a random variable. Then $\bH(X) \le \log_2  \lvert \supp(X) \rvert$, where the equality is attained if and only if the distribution of $X$ is uniform on $\supp(X)$. 
\end{proposition}

For any random vector $X = (X_1, \dots, X_n)$, we write $\bH(X_1, \dots, X_n) \eqdef \bH(X)$. For two random variables $X, Y$, we denote $p_X(x) \eqdef \p(X = x), \, p_Y(y) \eqdef \p(Y = y), \, p_{X, Y}(x, y) \eqdef \p(X = x, \, Y = y)$, respectively. Write $X \mid Y = y$ as the conditional random variable ``$X$ given $Y = y$'', and denote by $\bH(X \mid Y = y)$ its entropy, where $\bH(X \mid Y = y) = 0$ if $p_Y(y) = 0$. Define the conditional entropy
\[
\bH(X \mid Y) \eqdef \sum_{x, y} -p_{X, Y}(x, y) \log_2 \biggl( \frac{p_{X, Y}(x, y)}{p_Y(y)} \biggr) = \sum_{y} p_Y(y) \cdot \bH(X \mid Y = y). 
\]
One can easily verify that $\bH(X \mid Y) = \bH(X, Y) - \bH(Y)$. This implies the following ``chain rule'': 

\begin{proposition} \label{prop:entropy_chain}
    Let $X_1, \dots, X_n$ be random variables on the probability space $(\Omega, \p)$. Then
    \[
    \bH(X_1, \dots, X_n)=\bH(X_n \mid X_1, \dots, X_{n-1})+\bH(X_{n-1} \mid X_1, \dots, X_{n-2})+\dots+\bH(X_1). 
    \]
\end{proposition}

We will also need the following property which, informally speaking, says that adding already-known information does not increase the entropy.
\begin{proposition} \label{prop:entropy_deterministic_function}
    Let $X$ be a random variable on $(\Omega, \p)$. For any deterministic function $f$, we have 
    \[
    \bH(X) = \bH \bigl( X, f(X) \bigr). 
    \]
\end{proposition}

We introduce an inequality playing a central role in the entropic proof of the Kruskal--Katona theorem \cite{chao_yu}. For completeness, we also include its short proof. Write $[n] \eqdef \{1, 2, \dots, n\}$. 

\begin{lemma} \label{lem:entropy}
    Let $\mathcal{A}$ be a family of $d$-subsets of $[n]$. Uniformly at random, we sample a set $A \in \mathcal{A}$ and then uniformly randomly order the elements of $A$ into a vector $(X_1, \dots, X_d) \in [n]^d$. Then
    \[
    2^{\bH(X_1)} \ge 2^{\bH(X_2 \mid X_1)} + 1 \ge 2^{\bH(X_3 \mid X_1, X_2)} + 2 \ge \cdots \ge 2^{\bH(X_d \mid X_1, \dots, X_{d-1})} + d - 1. 
    \]
\end{lemma}

\begin{proof}
    It suffices to show $2^{\bH(X_k \mid X_1, \dots, X_{k-1})} \ge 2^{\bH(X_{k+1} \mid X_1, \dots, X_k)} + 1$ for all $k \in [d]$. Let $Z \in \{k, k+1\}$ be a random variable that is independent from any of $X_1, \dots, X_d$ such that
    \begin{align*}
        p \eqdef \bP(Z=k) &= \frac{2^{\bH(X_1, \dots, X_k)}}{2^{\bH(X_1, \dots, X_k)} + 2^{\bH(X_1, \dots, X_{k+1})}}, \\
        \bP(Z=k+1) &= \frac{2^{\bH(X_1, \dots, X_{k+1})}}{2^{\bH(X_1, \dots, X_k)} + 2^{\bH(X_1, \dots, X_{k+1})}}. 
    \end{align*}
    Since $Z$ is determined by revealing ``if $X_Z = X_k$ holds'', from \Cref{prop:entropy_deterministic_function} we deduce that
    \begin{align*}
        \bH(X_1, \dots, X_k, X_Z) &= \bH(X_1, \dots, X_k, X_Z, Z) \overset{(*)}{=} \bH(Z) + \bH(X_1, \dots, X_k, X_Z \mid Z) \\
        &= \bH \bigl( \Bern(p) \bigr) + p \bH(X_1, \dots, X_k, X_k) + (1-p) \bH(X_1, \dots, X_k, X_{k+1}) \\
        &= p \bigl( \bH(X_1, \dots, X_k) - \log_2 p \bigr) + (1-p) \bigl( \bH(X_1, \dots, X_{k+1}) - \log_2 (1-p) \bigr) \\
        &= \log_2 \bigl( 2^{\bH(X_1, \dots, X_k)} + 2^{\bH(X_1, \dots, X_{k+1})} \bigr), 
    \end{align*}
    where at the step marked with ($*$) we applied the chain rule (\Cref{prop:entropy_chain}). We thus obtain
    \begin{align*}
        2^{\bH(X_1, \dots, X_k)} \bigl( 2^{\bH(X_{k+1} \mid X_1, \dots, X_k)} + 1 \bigr) &= 2^{\bH(X_1, \dots, X_{k+1})} + 2^{\bH(X_1, \dots, X_k)} = 2^{\bH(X_1, \dots, X_k, X_Z)} \\
        &= 2^{\bH(X_1, \dots, X_k)} 2^{\bH(X_Z \mid X_1, \dots, X_k)} \le 2^{\bH(X_1, \dots, X_k)} 2^{\bH(X_Z \mid X_1, \dots, X_{k-1})}. 
    \end{align*}
    When conditioning on $X_1, \dots, X_{k-1}$, the distributions of $X_k$ and $X_Z$ are the same. So, 
    \[
    2^{\bH(X_{k+1} \mid X_1, \dots, X_k)} + 1 \le 2^{\bH(X_Z \mid X_1, \dots, X_{k-1})} = 2^{\bH(X_Z \mid X_1, \dots, X_{k-1})} = 2^{\bH(X_k \mid X_1, \dots, X_{k-1})}. \qedhere
    \]
\end{proof}

\section{The maximum number of \texorpdfstring{$t$}{t}-cliques for small \texorpdfstring{$p$}{p}} \label{sec:psmall}

We are ready to prove \Cref{thm:main}\ref{main_psmall}. The idea is to bound $\bH(X_1) + p\bH(X_2 \mid X_1)$ from both sides, as this quantity allows us to build up the connection between degree sequence $p$-norm and entropy. 
\vspace{-0.5em}
\begin{itemize}
    \item On one hand, we shall utilize \Cref{lem:entropy} to lower bound $\bH(X_1) + p\bH(X_2 \mid X_1)$. 
    \vspace{-0.5em}
    \item On the other hand, we shall upper bound $\bH(X_1) + p\bH(X_2 \mid X_1)$ in terms of $\|\bm{d}(G)\|_p$. 
\end{itemize}
\vspace{-0.5em}

\smallskip

If $0 \le u < t$, then $\binom{u}{t} < 1$, and \Cref{thm:main}\ref{main_psmall} trivially holds. We assume $u \ge t$ from now on and prove the contrapositive ``If $k_t(G) > \binom{u}{t}$, then $\|\bm{d}(G)\|_p > \sqrt[p]{u}(u-1) = C$.'' 

Let $A$ be a $t$-clique chosen uniformly at random and $(X_1, \dots, X_t)$ be a random ordering of the vertices of $A$ chosen uniformly. For $i \in [t]$, set $x_i \eqdef 2^{\bH(X_i \mid X_1, \dots, X_{i-1})}$. Since $(X_1, \dots, X_t)$ is uniform, it follows from the chain rule (\Cref{prop:entropy_chain}) and the uniform bound (\Cref{prop:entropy_max}) that 
\begin{equation} \label{eq:x_product}
    x_1 x_2 \cdots x_t = 2^{\bH(X_1, \dots, X_t)} = t! k_t(G) > u (u-1) \cdots (u-t+1). 
\end{equation}
Meanwhile, it follows from \Cref{lem:entropy} that
\begin{equation} \label{eq:x_monotone}
    x_1 \ge x_2 + 1 \ge \dots \ge x_t + t - 1. 
\end{equation}

\begin{claim} \label{claim:psmall_entropy_lower}
    For any $p \in [0, t-1]$, we have that $x_1 x_2^p > u(u-1)^p$. 
\end{claim}

\begin{proof}
    From \eqref{eq:x_monotone} we deduce the following estimates: 
    \begin{align}
        x_1^{} x_2^p &\ge (x_2^{}+1)x_2^p, \label{eq:x_entropy1} \\
        x_2^{1-p} x_3^{} \cdots x_t^{} &\le x_2^{1-p}(x_2^{}-1)\cdots(x_2^{}-t+2). \label{eq:x_entropy2}
    \end{align}
    Observe that both $(x_2^{}+1)x_2^p$ and $x_2^{1-p}(x_2^{}-1)\cdots(x_2^{}-t+2)$ are strictly increasing functions in $x_2$. The former is obvious, while the latter can be seen by rewriting the expression as 
    \[
    x_2^{t-1-p} \cdot \frac{x_2-1}{x_2} \cdot \dots \cdot \frac{x_2-t+2}{x_2}, 
    \]
    which is a product of increasing functions in $x_2$. Therefore, 
    \vspace{-0.5em}
    \begin{itemize}
        \item if $x_2 > u - 1$, then \eqref{eq:x_entropy1} implies that $x_1^{} x_2^p > u(u-1)^p$; 
        \vspace{-0.5em}
        \item if $x_2 \le u - 1$, then \eqref{eq:x_entropy2} and \eqref{eq:x_product} implies that $x_1^{} x_2^p > u(u-1)^p$. 
    \end{itemize}
    \vspace{-0.5em}
    Since the inequality holds in both cases, \Cref{claim:psmall_entropy_lower} is valid. 
\end{proof}

By taking logarithm of base $2$, \Cref{claim:psmall_entropy_lower} tells us that $\bH(X_1) + p\bH(X_2 \mid X_1) > \log_2 \bigl( u(u-1)^p \bigr)$. 

\begin{claim} \label{claim:entropy_upper}
    We have $\bH(X_1) + p\bH(X_2 \mid X_1) \le \log_2 \Bigl( \sum\limits_{v \in V} \deg(v)^p \Bigr)$. 
\end{claim}
\begin{proof}
    
Write $p_v \eqdef \bP(X_1 = v)$. \Cref{prop:entropy_max} implies that $\bH(X_2 \mid X_1 = v) \le \log_2 \bigl( \deg(v) \bigr)$, and hence
\begin{align*}
    \bH(X_1) + p\bH(X_2 \mid X_1) &\le \sum_{v \in V} \Bigl( -p_v \log_2 p_v + p_v \log_2 \bigl( \deg(v)^p \bigr) \Bigr) \\
    &= \sum_{v \in V} p_v \log_2 \biggl( \frac{\deg(v)^p}{p_v} \biggr) \\
    &\le \log_2 \biggl( \sum_{v \in V} \deg(v)^p \biggr), 
\end{align*}
where the second ``$\le$'' follows from the concavity of $f(x) = \log_2 x \, (x > 0)$.
\end{proof}

We thus conclude that
\[
\log_2 \biggl( \sum_{v \in V} \deg(v)^p \biggr) \ge \bH(X_1) + p\bH(X_2 \mid X_1) > \log_2 \bigl( u(u-1)^p \bigr), 
\]
which implies the desired inequality $\|\bm{d}(G)\|_p > \sqrt[p]{u}(u-1)$. The proof of \Cref{thm:main}\ref{main_psmall} is complete. 

\section{The maximum number of \texorpdfstring{$t$}{t}-cliques for large \texorpdfstring{$p$}{p}} \label{sec:plarge}

We now focus on the $p > t - 1$ case. Before proving the results, we first discuss some properties of 
\[
h(x) = \frac{1}{x(x-1)^p} \binom{x}{t}
\]
and the maximum points $s_{\R}, s_{\N}$. 

\begin{proposition} \label{prop:srn}
    If $t \ge 3$ and $p > t - 1$, then $s_{\R}$ is the unique critical point of $h(x)$ on $(t-1, +\infty)$. Moreover, $h(x)$ is strictly increasing on $(t-1, s_{\R}]$ and strictly decreasing on $[s_{\R}, +\infty)$. In particular, the maximum point $s_{\R}$ is well-defined, and there are at most two candidates for $s_{\N}$. 
\end{proposition}

\begin{proof}
    By taking the derivative, we have that
    \[
    h'(x) = \biggl( \sum_{i=1}^{t-1}\frac{1}{x-i}-\frac{p}{x-1} \biggr) h(x) = \biggl( \sum_{i=1}^{t-1}\frac{x-1}{x-i}-p \biggr) \frac{h(x)}{x-1} = \biggl( \sum_{i=1}^{t-1}\frac{i-1}{x-i}+t-1-p \biggr) \frac{h(x)}{x-1}. 
    \]
    Since $\frac{h(x)}{x-1} > 0$ holds for all $x > t-1$, it suffices to look at the sign of 
    \[
    g(x) = \sum_{i=1}^{t-1} \frac{i-1}{x-i} + t - 1 - p. 
    \]
    Note that $\lim\limits_{x \to (t-1)^+} g(x) = +\infty$ and $g(x) < 0$ if $x$ is large enough since $t - 1 - p < 0$. Moreover, it is easily seen that $g(x)$ is strictly decreasing. So, there exists a unique critical point $s_{\R}$ satisfying $g(s_{\R}) = 0$ that maximizes the function $h(x)$. We thus conclude the proposition. 
\end{proof}

\subsection{Proof of \texorpdfstring{\Cref{thm:clique_plarge}}{Theorem 4}}

Again, we prove the contrapositive ``If $k_t(G)>\frac{n}{u}\binom{u}{t}$, then $\|\bm{d}(G)\|_p > \sqrt[p]{n}(u-1) = C$.'' The idea is similar to that behind \Cref{sec:psmall}. Sample $X_1, \dots, X_t$ and define $x_1, \dots, x_t$ in the same way we did in \Cref{sec:psmall}. Due to the chain rule (\Cref{prop:entropy_chain}) and the uniform bound (\Cref{prop:entropy_max}), we have
\begin{equation} \label{eq:clique_plarge_x_product}
    x_1 x_2 \cdots x_t = 2^{\bH(X_1, \dots, X_t)} = t! k_t(G) > n (u-1)(u-2) \cdots (u-t+1). 
\end{equation}
Meanwhile, it follows from \Cref{lem:entropy} and the uniform bound (\Cref{prop:entropy_max}) that
\begin{equation} \label{eq:clique_plarge_x_monotone}
    n\ge x_1 \ge x_2 + 1 \ge \dots \ge x_t + t - 1. 
\end{equation}
Also, from $n(s_{\R}-1)^p \le C^p$ we deduce that $s_{\R} \le \frac{C}{\sqrt[p]{n}} + 1 = u$. 

\begin{claim} \label{claim:nplarge_entropy_lower}
    For any $p > t - 1$, we have $x_1 x_2^p > n(u-1)^p$. 
\end{claim}

\begin{proof}
By combining \eqref{eq:clique_plarge_x_product} and \eqref{eq:clique_plarge_x_monotone}, we obtain
\[
x_1x_2(x_2-1)\cdots(x_2-t+2) > n(u-1)(u-2) \cdots (u-t+1). 
\]
By rearranging the above inequality, the fact $x_1 \le n$ implies
\[
1 \ge \frac{x_1}{n} > \frac{(u-1)(u-2) \cdots (u-t+1)}{x_2(x_2-1) \cdots (x_2-t+2)}, 
\]
and so $s_{\R} \le u < x_2+1$. \Cref{prop:srn} tells us that $h(x)$ is strictly decreasing on $[s_{\R}, +\infty)$, and so
\[
\frac{1}{u(u-1)^p} \binom{u}{t} = h(u) > h(x_2+1) = \frac{1}{(x_2+1)x_2^p} \binom{x_2+1}{t}. 
\]
It follows that
\[
\frac{x_1}{n} > \frac{(u-1)(u-2) \cdots (u-t+1)}{x_2(x_2-1) \cdots (x_2-t+2)} > \frac{(u-1)^p}{x_2^p} \implies x_1x_2^p> n(u-1)^p. \qedhere
\]

By taking logarithm of base $2$, \Cref{claim:nplarge_entropy_lower} tells us that 
\[
\bH(X_1) + p\bH(X_2 \mid X_1) > \log_2 \bigl( n(u-1)^p \bigr). 
\]
Since \Cref{claim:entropy_upper} still holds here, it follows that
\[
\log_2 \biggl( \sum_{v \in V} \deg(v)^p \biggr) \ge \bH(X_1) + p\bH(X_2 \mid X_1) > \log_2 \bigl( n(u-1)^p \bigr), 
\]
which implies the desired inequality $\|\bm{d}(G)\|_p > \sqrt[p]{n}(u-1)$. The proof of \Cref{thm:clique_plarge} is complete. 
\end{proof}

\subsection{Proof of \texorpdfstring{\Cref{thm:main}\ref{main_plarge}}{Theorem 3(ii)}}

Let the vertices of $G$ be $v_1, \dots, v_n$. Set $a_i \eqdef \deg(v_i)$. A direct double counting shows that 
\begin{align*}
    k_t(G) &\le \frac{1}{t} \sum_{i=1}^n \binom{\deg(v_i)}{t-1} = \frac{1}{t} \sum_{i=1}^n \binom{a_i}{t-1}= \sum_{i=1}^n \frac{1}{a_i+1} \binom{a_i+1}{t} \\
    &= \sum_{i=1}^n h(a_i+1) a_i^p \le h(s_{\N}) \sum_{i=1}^n a_i^p \le h(s_{\N}) C^p = \frac{C^p}{u(u-1)^p} \binom{u}{t}. 
\end{align*}

\section{Concluding remarks} \label{sec:remark}

Our results can be extended to hypergraphs, which generalizes the original Kruskal--Katona Theorem concerning bounding the number of $K^{(r)}_t$-hypercliques in an $r$-uniform hypergraph. 

Suppose $p \in (0, \infty)$ and $r > j$ are positive integers. Let $H = (V, E)$ be an $r$-uniform hypergraph. For any $j$-subset $S \in \binom{V}{j}$, we denote by $\deg(S)$ the number of $(r-j)$-subsets $T \in \binom{V \setminus S}{r-j}$ satisfying $S \cup T \in E$. This generalizes the usual degree of a vertex in a graph. Define the $(j, p)$-norm of $H$ as
\[
\|\bm{d}(H)\|_{j, p} = \Biggl(\sum_{S \in \binom{V}{j}} \deg(S)^p \Biggr)^{\frac{1}{p}}. 
\]
This definition is essentially the same as the $(t, p)$-norm studied in \cite{chen_ilkovic_leon_liu_pikhurko} (with $j$ in the place of $t$). 

We are to maximize the number of $r$-uniform $t$-clique subgraphs in $H$, denoted by $k^r_t(H)$, given that the $(j, p)$-norm of $H$ is bounded from the above. Here the results and the proofs for general $r, j$ with $t > r > j$ when $p$ is small are parallel to those for the $(r, j) = (2, 1)$ special case, i.e., \Cref{thm:main}\ref{main_psmall}. Define
\[
\widetilde{h}(x) = \frac{\binom{x}{t}}{\binom{x}{j}\binom{x-j}{r-j}^p}
\]
on the interval $(t-1, +\infty)$. When $p > \frac{t-j}{r-j}$, set $\widetilde{s}_{\R}$ to be the real number maximizing $h(x)$. Similar to \Cref{prop:srn}, a derivative argument shows the existence and the uniqueness of $\widetilde{s}_{\R}$. Moreover, the function $\widetilde{h}(x)$ is strictly increasing on $(t-1, \widetilde{s}_{\R}]$ and strictly decreasing on $[\widetilde{s}_{\R}, +\infty)$. 

\begin{theorem} \label{thm:main_hyper}
    Suppose $C > 0$. Let $H$ be an $r$-uniform hypergraph with $\|\bm{d}(H)\|_{j, p} \le C$. 
    \vspace{-0.25em}
    \begin{enumerate}[label=(\roman*), ref=(\roman*)]
        \item \label{hyper_psmall} If $0 < p \le \frac{t-j}{r-j}$, then $k_t^r(H) \le \max \bigl\{ \binom{u}{t}, 0 \bigr\}$, where $u$ is the positive real with $C = \sqrt[p]{\binom{u}{j}} \binom{u-j}{r-j}$. 
        \vspace{-0.5em}
        \item \label{hyper_plarge} If $p > \frac{t-j}{r-j}$, then $k_t^r(H) \le \frac{C^p\binom{u}{t}}{\binom{u}{j}\binom{u-j}{r-j}^p}$, where $u = \widetilde{s}_{\R}$. 
    \end{enumerate}
    \vspace{-0.25em}
\end{theorem}

\Cref{thm:main_hyper} is tight, as witnessed by examples similar to those of \Cref{thm:main}. The proof idea is the same as that behind \Cref{thm:main}. However, the detailed estimates inside turn out to be more technical, and we include a sketched proof in \Cref{append:main_hyper}. It is worth mentioning that studying the tight upper bound on $k_t^r(H)$ when the number of vertices is fixed could also be interesting. 

\smallskip

Recall that \Cref{thm:GLS_chase} answers \Cref{ques:triangle_maxdeg} under the imposed assumption that $G$ is an $n$-vertex graph (vertex-problem). Kirsch and Radcliffe \cite{kirsch_radcliffe_2019} proposed a natural variant problem where they fix the number of edges (edge-problem) rather than the number of vertices. This problem was completely solved by Chakraborti and Chen \cite{chakraborti_chen}. Later, Kirsch and Radcliffe \cite{kirsch_radcliffe_2023} proved generalizations of both the vertex-problem and the edge-problem in hypergraphs. For future work, it could be an interesting problem to maximize $k_t(G)$, given that $\|\bm{d}(G)\|_p$ and $\|\bm{d}(G)\|_q$ are fixed for some $p < q$. Notice that the aforementioned edge-problem is a special case of this framework, where $p = 1$ and $q = \infty$. 

\section*{Acknowledgments}

We are grateful to Hong Liu for carefully reading a draft of this paper and providing helpful suggestions. We also benefited from discussions with Debsoumya Chakraborti and Ruonan Li. 

This work was initiated at the $2^{\text{nd}}$ ECOPRO Student Research Program in the summer of 2024. The third and the fourth authors would like to thank ECOPRO for hosting them. 

\bibliographystyle{plain}
\bibliography{norm_clique}

\appendix

\section{Proof of \texorpdfstring{\Cref{thm:main_hyper}}{Theorem 10}} \label{append:main_hyper}

Randomly pick and label a $t$-clique into $(X_1, \dots, X_t)$ as before. One can prove the estimates
\begin{align}
    \log_2 \Bigl( j!(r-j)!^p \sum_{S \in \binom{V}{j}} \deg(S)^p \Bigr) &\ge \bH(X_1, \dots, X_j) + p\bH(X_{j+1}, \dots, X_r \mid X_1, \dots, X_j) \label{eq:hyperentropy1} \\
    &> \log_2 \Bigl( \bigl( u(u-1)\cdots(u-j+1) \bigr) \bigl( (u-j)\cdots(u-r+1) \bigr)^p \Bigr). \label{eq:hyperentropy2}
\end{align}
Evidently, this chain of inequalities will conclude \Cref{thm:main_hyper}. 

Similar to the proof of \Cref{thm:main}\ref{main_psmall}, the inequality \eqref{eq:hyperentropy1} follows from the observation that
\[
\bH(X_{j+1}, \dots, X_r \mid X_1 = v_1, \dots, X_j = v_j) \le \log_2 \Bigl( (r-j)!\deg \bigl( \{v_1, \dots, v_j\} \bigr) \Bigr). 
\]
To see the estimate \eqref{eq:hyperentropy2}, we write $x_i \eqdef 2^{\bH(X_i \mid X_1, \dots, X_{i-1})}$ and deduce the following as before: 
\begin{equation*}
    x_1 x_2 \cdots x_t = 2^{\bH(X_1, \dots, X_t)} = t! k_t^r(H) > u (u-1) \cdots (u-t+1). 
\end{equation*}
\begin{equation*}
    x_1 \ge x_2 + 1 \ge \dots \ge x_t + t - 1. 
\end{equation*}

Write $A \eqdef x_1 \cdots x_j, \, B \eqdef x_{j+1} \cdots x_r, \, C \eqdef x_{r+1} \cdots x_t$. It then suffices to prove that
\[
AB^p \ge \bigl( u(u-1) \cdots (u-j+1) \bigr) \bigl( (u-j)(u-j-1) \cdots (u-r+1) \bigr)^p. 
\]
If $B \ge (u-j) \cdots (u-r+1)$, then $x_{j+1} \ge u-j$, and so $x_1 \ge u, \dots, x_j \ge u-j+1$. It follows that
\[
AB^p \ge \bigl( u \cdots (u-j+1) \bigr) \bigl( (u-j) \cdots (u-r+1) \bigr)^p. 
\]
If $B \le (u-j) \cdots (u-r+1)$, then $x_r \le u-r+1$. So, 
\[
AB^p = \frac{ABC}{B^{1-p}C} = \frac{u \cdots (u-t+1)}{B^{1-p}C}. 
\]
Find $b \ge t$ such that $B = (b-j)\cdots(b-r+1)$. Then $u \ge b$. Observe that
\[
B = x_{j+1} \cdots x_r \ge (x_r+r-j-1) \cdots x_r. 
\]
It follows that $b \ge x_r + r - 1$. Therefore, from $x_{r+1} \cdots x_t \le (b-r) \cdots (b-t+1)$ we deduce that
\[
B^{1-p}C \le \bigl( (b-j) \cdots (b-r+1) \bigr)^{1-p} (b-r)\cdots(b-t+1). 
\]
Rewrite the right-hand side as
\[
\bigl( (b-j) \cdots (b-r+1) \bigr)^{\frac{t-r}{r-j}+1-p} \frac{b-r}{\bigl( (b-j) \cdots (b-r+1) \bigr)^{\frac{1}{r-j}}} \cdots \frac{b-t+1}{\bigl( (b-j) \cdots (b-r+1) \bigr)^{\frac{1}{r-j}}}. 
\]
\begin{itemize}
    \item The first part is increasing in $b$ because $p \le \frac{t-j}{r-j}$, and 
    \vspace{-0.5em}
    \item the second part is increasing in $b$ since $f(x) = \frac{x-\lambda}{x-\mu} \, (x > \lambda)$ is increasing whenever $\lambda > \mu$. 
\end{itemize}
\vspace{-0.5em}
We obtain $B^{1-p}C \le \bigl( (u-j) \cdots (u-r+1) \bigr)^{1-p} (u-r) \cdots (u-t+1)$, and hence conclude that
\[
AB^p = \frac{u \cdots (u-t+1)}{B^{1-p}C} \ge \bigl( u \cdots (u-j+1) \bigr) \bigl( (u-j) \cdots (u-r+1) \bigr)^p. 
\]

The proof of \Cref{thm:main_hyper}\ref{hyper_psmall} is complete, and we focus on \Cref{thm:main_hyper}\ref{hyper_plarge} then. 

For any $S \in \binom{V}{j}$, let $a_S\geq r-j$ be the positive real with $\deg(S) = \binom{a_S}{r-j}$. The Lov\'{a}sz approximate version of the Kruskal--Katona theorem tells us that, by restricting $H$ to obtain an $\binom{a_S}{r-j}$-hyperedge $(r-j)$-uniform link hypergraph on the neighborhood of $S$, there are at most $\binom{a_S}{t-j}$ many $(t-j)$-hypercliques. A direct double counting argument then shows that
\begin{align*}
    k_t^r(H) &\le \frac{1}{\binom{t}{j}} \sum_{S \in \binom{V}{j}} \binom{a_S}{t-j} = \sum_{S \in \binom{V}{j}} \frac{\binom{a_S+j}{t}}{\binom{a_S+j}{j}} = \sum_{S \in \binom{V}{j}} \widetilde{h}(a_S+j) \binom{a_S}{r-j}^p \\
    &\le \widetilde{h}(\widetilde{s}_{\R}) \sum_{S \in \binom{V}{j}} \binom{a_S}{r-j}^p \le \widetilde{h}(\widetilde{s}_{\R}) C^p = \frac{C^p\binom{u}{t}}{\binom{u}{j}\binom{u-j}{r-j}^p}, 
\end{align*}
which concludes the proof of \Cref{thm:main_hyper}\ref{hyper_plarge}. 

\end{document}